\documentclass[
final
]{dmtcs-episciences}

\usepackage[utf8]{inputenc}
\usepackage{subfigure}
\usepackage{enumerate}
\usepackage{amsthm}
\usepackage{amsmath}

\newtheorem{theorem}{Theorem}[section]
\newtheorem{conjecture}[theorem]{Conjecture}
\newtheorem*{burning}{Burning number conjecture}
\newtheorem{question}[theorem]{Question}
\newtheorem{observation}[theorem]{Observation}
\theoremstyle{definition}
\newtheorem{definition}[theorem]{Definition}

\theoremstyle{remark}
\newtheorem{example}[theorem]{Example}
\newtheorem{remark}[theorem]{Remark}

\usepackage[round]{natbib}

\author[Ta Sheng Tan and Wen Chean Teh]{Ta Sheng Tan\affiliationmark{1}
  \and Wen Chean Teh\affiliationmark{2}\thanks{Corresponding author} }
\title[A note on graph burning of path forests]{
A note on graph burning of path forests
}
\affiliation{
  Institute of Mathematical Sciences, Faculty of Science,  University Malaya,  Malaysia\\
  School of Mathematical Sciences,  Universiti Sains Malaysia,  Malaysia}
\keywords{discrete graph algorithm, burning number conjecture, spread of social contagion, sumset partition of integers, well-burnable}
\begin{document}
\publicationdata{vol. 26:3}{2024}{1}{10.46298/dmtcs.12709}{2023-12-19; 2023-12-19; 2024-03-21}{2024-07-12}

\maketitle
\begin{abstract}
 Graph burning is a natural discrete graph algorithm inspired by the spread of social contagion. Despite its simplicity, some open problems remain steadfastly unsolved, notably the burning number conjecture, which says that every connected graph of order $m^2$ has burning number at most $m$. Earlier, we showed that the conjecture also holds for a path forest, which is disconnected, provided each of its paths is sufficiently long. However, finding the least sufficient length for this to hold turns out to be nontrivial. In this note, we present our initial findings and conjectures that associate the problem to some naturally impossibly burnable path forests. It is noteworthy that our problem can be reformulated as a topic concerning sumset partition of integers.
\end{abstract}

\section{Introduction}\label{intro}

Graph burning is a discrete-time process introduced by \citet{bonato2016how}
that can be viewed as a simplified model for the spread of contagion in a network.
Given a simple finite graph $G$, each vertex of the graph is either \emph{burned} or \emph{unburned} throughout the process.
Initially, every vertex of $G$ is unburned.
At the beginning of every round $t\ge 1$, a \emph{burning source} is placed at an unburned vertex to burn it.
If a vertex is burned in round $t-1$, then in round $t$, each of its unburned neighbours becomes burned.
A burned vertex will remain burned throughout the process.
The burning process ends when all vertices of $G$ are burned, in which case we say the graph $G$ is \emph{burned}.
The \emph{burning number of $G$} is the least number of rounds needed for the burning process to be completed.

The study of graph burning is extensive, with the main open problem being the burning number conjecture by~\citet{bonato2016how}.

\begin{burning}\citep{bonato2016how}
	The burning number of every connected graph of order $N$ is at most $\lceil\sqrt{N}\rceil$.
\end{burning}

In the literature of graph burning, a graph is said to be \emph{$m$-burnable} if its burning number is at most $m$, and a graph (including a disconnected graph) is said to be \emph{well-burnable} if it satisfies the burning number conjecture.
Many classes of graphs have been verified to be well-burnable, including hamiltonian graphs \citep{bonato2016how}, spiders \citep{bonato2019bounds, das2018burning}, and caterpillars \citep{hiller2021burning, liu2020burning}. Recently, the burning number conjecture was shown to hold asymptotically by \cite{norin2024burning}.

The reader is referred to \citet{bonato2021survey} for a survey on graph burning. In this short note, we are interested in the graph burning of path forests.
Here, a \emph{path forest} is a disjoint union of paths.
While not all path forests are well-burnable, it was shown in \citet{tan2023burnability} that a path forest with a sufficiently long shortest path is well-burnable.

\begin{theorem}\citep{tan2023burnability}\label{030922c}
	For every $n\in \naturals$, there exists a smallest $L_n\in \naturals$ such that
	if $T$ is a path forest with $n$ paths and the shortest path of $T$ has order at least $L_n$, then $T$ is well-burnable.
\end{theorem}

When determining whether a path forest is well-burnable, we can always extend some of its paths so that the order of the graph is $m^2$ for some $m\in \naturals$.
So from here onwards, we will assume that the order of a path forest is always an integer squared.
We say a path forest is \emph{deficient} if it is not well-burnable, and by an \emph{$n$-path forest}, we mean a path forest with $n$ paths.

The main purpose of this note is to provide insights on $L_n$, and some conjectures related to $L_n$.
We clearly have $L_1 = 1$ as every path is well-burnable, and it is straightforward that $L_2 = 3$, since a $2$-path forest is deficient if and only if its path orders are $m^2-2$ and $2$ for some $m\geq 2$ \citep[see][Lemma~3.1]{tan2020graph}.
The study of the values of $L_n$ was posed as an open problem in \citet{tan2023burnability}, where it was mentioned that $L_3=18$ and $L_4=26$ (determined with careful analysis and the help of a computer).

Note that for a path forest of order $m^2$ to be well-burnable, a burning process of $m$ rounds must have its $i$th burning source burning exactly $2m-2i+1$ vertices, and every vertex is burned by exactly one burning source.
So, each path in the path forest is burned exactly, and the number of burning sources used on the path has the same parity as its order.
While investigating deficient $n$-path forests for some small values of $n$, we notice that when their shortest paths have orders slightly smaller than $L_n$, they are all trivially deficient, in the following sense.
For such a deficient path forest, even if we pretend that for each of its paths, the $i$th burning source used on it (not on the path forest) would burn $2m-2i+1$ vertices, $m$ burning sources are not enough to completely burn the path forest, provided we insist that the number of burning sources used on each path has the same parity as its order.
We will call these trivially deficient path forests \emph{impossibly burnable} (this will be made more precise in Section 2).

\begin{question}\label{mainQuestion}
	Is it true that if $T$ is a deficient $n$-path forest such that its shortest path has order $L_n-1$, then $T$ is impossibly burnable?
\end{question}

As mentioned earlier, the answer to the above question is affirmative for small values of $n$ (up to $n=7$).
However, what about the remaining values of $n$?
While we are unable to verify for any $n\ge 8$ (due to computational limitations), we believe that the answer to the above question remains affirmative (see Conjecture~\ref{050522f}).
This belief has motivated us to study impossibly burnable path forests, which would lead to the computation of the exact values of $L_n$.

For $n\ge 2$, let $M_n$ be the smallest positive integer such that if $T$ is an impossibly burnable $n$-path forest, then its shortest path has order at most $M_n$.
So assuming the answer to Question~\ref{mainQuestion} is affirmative, it follows that $L_n = M_n + 1$ for all $n\ge 2$.
It is straightforward that $M_2 = 2$, and for the main result of this note, we determine the exact values of $M_n$.

\begin{theorem}\label{mainTheorem}
	For $n\ge 3$, $M_n$ is the largest odd number smaller than or equal to $$ 12n - 2\sqrt{18n - 12} - 6.$$
\end{theorem}

We digress briefly to mention another formulation of the graph burning of path forests, presented in the form of a sumset partition problem.
The problem of determining whether an $n$-path forest $(l_1,l_2,\ldots,l_n)$ of order $m^2$ is well-burnable is equivalent to deciding whether the set of the first $m$ odd positive integers can be partitioned into $n$ subsets $S_1, S_2, \ldots, S_n$ such that for every $i\in [n]$, the sum of the numbers in $S_i$ is equal to $l_i$.
Interested readers may refer to \citet{ando1990disjoint}, \citet{chen2005partition}, \citet{ enomoto1995disjoint}, \citet{fu1994disjoint}, \citet{llado2012modular}, and \citet{ma1994ascending} for some related studies on this formulation.

\section{Impossibly burnable path forests}

In this section, we will first describe our observations on deficient path forests that lead us to Question~\ref{mainQuestion}, and then we will proceed to prove Theorem~\ref{mainTheorem}.
For a path forest, its \emph{path orders} indicate the respective order of each of its paths.
We may represent an $n$-path forest by an $n$-tuple $(l_1,l_2,\cdots,l_n)$ of its path orders.
Often, we assume $l_1\le l_2\le \cdots \le l_n$, and if $T$ is a path forest, we may write $l_1(T)$ (or just $l_1$) for the order of its shortest path.

In \cite{bessy2017burning}, the graph burning problem was shown to be NP-complete for general path forests, and a polynomial time algorithm for the problem was constructed when the number of paths is fixed. Based on this algorithm, given the number of paths $n$ and a positive integer $M\geq n$, we are able to construct  a complete list of well-burnable $n$-path forests of order $m^2$ for all $m\leq M$ using Matlab subject to computational limitations. Note that $(1,3,5, \ldots, 2n-1)$ is  the unique well-burnable $n$-path forest of order $n^2$. Our lists are constructed recursively based on the following strategy: to obtain a well-burnable $n$-path forest of order $(m+1)^2$,
\begin{enumerate}
	\item add a path of order $2m+1$ to any well-burnable $(n-1)$-path forest $T$ such that $\vert T\vert=m^2$; or
	\item add $2m+1$ vertices to any one of the paths of any well-burnable $n$-path forest $T$ such that $\vert T\vert=m^2$.
\end{enumerate}

Suppose $T$ is a $3$-path forest of order $m^2$ and $l_1\geq 8$.
From the complete list of well-burnable $3$-path forests for $m\leq 9$, we observe that if $T$ is deficient, then $T$ is one of the following six possibilities.
Furthermore, if $m=9$, then $T$ is well-burnable.
We can then deduce by induction that $T$ is well-burnable for any $m\geq 9$ by considering the $3$-path forest $(l_1,l_2,l_3-2m+1)$.

\begin{observation}\label{261021b}
	Every $3$-path forest with $l_1\geq 8$ is well-burnable, unless it is one of the following exceptional cases:
	$$(8,13,15), (8,15,26), (10,13,13), (15,15,19), (15, 17, 17), (17, 17, 30).$$
\end{observation}

Similarly, we have the following observation for $4$-path forests.

\begin{observation}\label{230222b}
	Every $4$-path forest with $l_1\geq 25$ is well-burnable, unless it is one of the following exceptional cases:
	$$(25,25,25,25), (25,25,25,46), (25,25,27,44), (25,25, 29, 42), (25,27,27,42), (25,25,25,69),$$
	$$ (25,25,27,67), (25,25, 29, 65), (25,27,27,65), (25,25,46,48), (25,27,46,46).$$	
\end{observation}

It follows from Observations~\ref{261021b} and \ref{230222b} that $L_3=18$ and $L_4=26$. Furthermore, as mentioned in the \nameref{intro}, we observed that all the deficient path forests in Observations~\ref{261021b} and \ref{230222b} are impossibly burnable.
We now give the precise definition of an impossibly burnable path forest.

\begin{definition}
	For $m\in\naturals$ and $1\le l\le m^2$, let $B_m(l)$ be the least $t\in\naturals$ having the same parity as $l$ such that $l\leq \sum_{i=1}^t [2m - (2i-1)] = 2mt-t^2$.
\end{definition}

\begin{example}
	If $l$ is odd and $ 2m-1< l \leq (2m-1)+(2m-3)+(2m-5)$, then $B_m(l)=3$.
	Meanwhile, if $l$ is even and $(2m-1)+(2m-3)< l \leq  (2m-1)+(2m-3)+(2m-5)+(2m-7)$, then $B_m(l)=4$.
\end{example}

\begin{definition}
	Suppose $T=(l_1,l_2,\ldots,l_n)$ is an $n$-path forest of order $m^2$.
	We say that $T$ is \emph{impossibly burnable} if $\sum_{i=1}^n B_m(l_i)>m$.
\end{definition}

\begin{remark}\label{040222a}
	\begin{enumerate}[(i)]
		\item An impossibly burnable path forest is clearly deficient. Indeed, for a path of order $l$ in a path forest of order $m^2$, at least $B_m(l)$ burning sources are required to burn the path completely in $m$ rounds, regardless of the parity of $l$.
		\item Not all deficient path forests are impossibly burnable. For example, the path forest $(2,7,7)$ is deficient but not impossibly burnable.
		\item The parities of $\sum_{i=1}^n B_m(l_i)$ and $m$ are equal, and thus if $T$ is impossibly burnable, then $m \leq \sum_{i=1}^n B_m(l_i)-2$.
	\end{enumerate}
\end{remark}

Based on Observation~\ref{261021b} and the subsequent discussion, any deficient $3$-path forest with $l_1\ge 8$ is impossibly burnable.
For the case of $4$-path forests, a Matlab search reveals that there are exactly $47$ deficient $4$-path forests with $l_1\geq 18$, all of which are impossibly burnable.
Note, however, that the path forest $(17,17,17,30)$ is deficient, but it is not impossibly burnable.
Analysing the cases for $5\le n\le 7$ gives similar results, leading us to the following conjecture.

\begin{conjecture}\label{050522f}
	Let $n\geq 4$.
	If $T$ is a deficient $n$-path forest with $l_1\geq L_{n-1}$, then $T$ is impossibly burnable.
\end{conjecture}

As mentioned earlier, Conjecture~\ref{050522f} is true for $n=4$, but it is not true for $n=3$ as $L_2=3$. The path forest $(3,3,3)$ is deficient but not impossibly burnable. Through an extensive Matlab search, we have determined that $L_5=36$, $L_6=46$, and $L_7=56$. Here, we briefly mention our computational validation of Conjecture~\ref{050522f}.

Analysing the list of well-burnable $5$-path forests order $m^2$ for $m\le 18$, we observe that there are exactly $608$ deficient $5$-path forests with $l_1\ge 26$, all of which are impossibly burnable.
Furthermore, all $5$-path forests of order $18^2$ with $l_1\ge 26$ are well burnable.
Hence, we can again deduce by induction that for $m\ge 18$, every $5$-path forest of order $m^2$ with $l_1\ge 26$ is well-burnable.
Similarly, Conjecture~\ref{050522f} holds true for $n=6$.
Specifically, all the $5185$ deficient $6$-path forests with $l_1\ge 36$ are impossibly burnable.
We have also managed to verify Conjecture~\ref{050522f} for $n=7$, with a more significant effort due to computational limitations.
(See \nameref{appendix A} for a brief account of this verification.)

As $M_n$ increases as $n$ grows by Theorem~\ref{mainTheorem}, we remark that Conjecture~\ref{050522f} implies a strongly affirmative answer to Question~\ref{mainQuestion}, resulting in $L_n = M_n+1$. We are now ready to determine the exact values of $M_n$.

\begin{proof}[of Theorem~\ref{mainTheorem}]
	Suppose $T$ is an impossibly burnable $n$-path forest of order $m^2$ with path orders $l_1\le l_2\le \cdots \le l_n$.
	Writing $t_i=B_m(l_i)$ for each $i\in [n]$ and recalling that $t_i\equiv l_i\pmod 2$ for every $i\in [n]$, we have that $\sum_{i=1}^n t_i \ge m+2$.
	Consider the partition of $[n]$ into
	\begin{displaymath}
	A = \{i\in[n]:t_i\ge 4\} \mbox{ and } B = \{i\in [n]:t_i\le 3\}.
	\end{displaymath}
	For convenience, we let $s_i = t_i - 2\ge 2$ for each $i\in A$.
	So for every $i\in A$, we have
	\begin{displaymath}
	l_i  \ge  (2m-1) + (2m-3) + \cdots + (2m - 2s_i + 1) + 2 =  2ms_i - s_i^2 + 2.
	\end{displaymath}
	
	Let $s = \sum_{i\in A} s_i = \left(\sum_{i\in A}t_i\right) - 2|A|$ and note that $s<m$, as otherwise, $\sum_{i\in A}l_i > m^2$.
	Observe now that
	\begin{eqnarray*}
	\sum_{i\in A} l_i &\ge & \sum_{i\in A} (2ms_i - s_i^2 +2)\\
	& = & 2ms - \sum_{i\in A} s_i^2 + 2|A|\\
&	=  & 2ms - s^2 + \left(\sum_{i,j\in A, i\ne j}s_is_j\right)+2|A|.
	\end{eqnarray*}
	It follows that $m^2 = \sum_{i=1}^n l_i \ge l_1|B| + 2ms - s^2 + \left(\sum_{i,j\in A,i\ne j }s_is_j\right)+2|A|$, implying that
	\begin{displaymath}
	(m-s)^2 \ge l_1|B| + \left(\sum_{i,j\in A, i\ne j}s_is_j\right) + 2|A|.
	\end{displaymath}
	
	On the other hand, note that $m+2 \le \sum_{i=1}^n t_i \le 3|B| + s + 2|A|$, or in other words,
	\begin{displaymath}
	0 < m - s \le 3|B| + 2|A| - 2.
	\end{displaymath}
	Putting these two inequalities together, we get
	\begin{equation}\label{mainInq}
	\left(3|B| + 2|A| - 2\right)^2 \ge l_1|B| + \left(\sum_{i,j\in A, i\ne j}s_is_j\right) + 2|A|.
	\end{equation}
	
	To bound $l_1$ from above, we consider a few cases. If $|B|=0$, we have \begin{math}(2n-2)^2 \ge \sum_{i,j\in [n], i\ne j} s_is_j + 2n \ge 4n(n-1) + 2n,\end{math} which is impossible, and so we must have $|B|>0$.
	If $|A| = 0$, we have $(3n-2)^2 \ge nl_1$, and so $l_1 \le 9n - 12 + \frac{4}{n}$.
	If $|A| = 1$, we have $(3n-3)^2 \ge (n-1)l_1 + 2$, and so $l_1 \le 9n - 9 - \frac{2}{n-1}$.
	
	For the final case where $|A|\ge 2$, we first observe that
	\begin{eqnarray*}
	\sum_{i,j\in A, i \ne j} s_is_j &= & \sum_{i\in A}s_i\left(\sum_{j\in A, j\ne i} s_j\right)
	= \sum_{i\in A} s_i (s-s_i)   \\
	&\ge & \sum_{i\in A} 2(s-2) \qquad\mbox{(as $2\le s_i\le s-2$ for every $i\in A$)}  \\
	&= &2|A|(s-2),
	\end{eqnarray*}
	Letting $s = 2|A| + k$ for some $k\ge 0$, we see that
	\begin{displaymath}
\sum_{i,j\in A, i \ne j} s_is_j\ge 2|A|(2|A| + k -2) = 4|A|^2+2k|A| - 4|A|.
\end{displaymath}
	Together with Inequality~(\ref{mainInq}), we have
	\begin{eqnarray*}
	(3|B|+2|A|-2)^2 &\ge & l_1|B| + 4|A|^2+2k|A| - 2|A|\\
	9|B|^2 + 12|A||B| - 12|B| - 6|A| + 4 &\ge & l_1|B| + 2k|A|\\
	\Longrightarrow\qquad \frac{2k|A|}{|B|} + l_1 &\le & 9|B|+12|A| - 6 - \frac{6(|A|+|B|) - 4}{|B|}\\
	\Longrightarrow\qquad \qquad \;\;\;\quad l_1 &\le &  9n -6 +3|A| - \frac{6n-4}{n-|A|}.
	\end{eqnarray*}
	It is straightforward that in the range of $0<x<n$, the function $3x - \frac{6n-4}{n-x}$ is maximised when $x = n - \sqrt{\frac{6n-4}{3}}$, with the maximum value being $3n - 2\sqrt{18n-12}$. Therefore, $l_1\le 12n - 2\sqrt{18n - 12} - 6$.
	
	We now see that for $n\ge 3$,
	\begin{eqnarray*}
	l_1 &\le & \max \left\{9n-12+\frac{4}{n}, 9n-9-\frac{2}{n-1},  12n - 2\sqrt{18n - 12} - 6\right\}  \\
&	= &12n - 2\sqrt{18n - 12} - 6.
	\end{eqnarray*}
	
	Before we proceed, we make a relevant observation. Pick 
	\begin{displaymath}	
	x_0\in \left\{\left\lfloor n-\sqrt{\frac{6n-4}{3}}\right\rfloor, \left\lceil n -\sqrt{\frac{6n-4}{3}}\,\right\rceil\right\}
	\end{displaymath}
	such that $3x - \frac{6n-4}{n-x}$ attains the larger value. 	It can be verified carefully but elementarily that the largest odd integer smaller than or equal to $9n -6 +3x_0 - \frac{6n-4}{n-x_0}$ coincides with that of $12n - 2\sqrt{18n - 12} - 6$ for every $n\geq 3$.

	Now, consider the $n$-path forest $T'$ with $m=3n+x_0-2$ and path orders as follows: 
	\begin{enumerate}
		\item 	$l'_i = 4m-2$ (and so $B_m(l_i')=4$) for each $i>n-x_0$;
		\item $l'_1, l'_2, \dots, l'_{n-x_0}$ are odd and any two of them are equal or differ by two. 	
	\end{enumerate} 	
	Such a path forest exists as the second requirement can be satisfied because $m$ is odd if and only if $n-x_0$ is odd. Note that $\sum_{i=1}^{n-x_0} l'_i$ is equal to
	\begin{eqnarray*}
	m^2 - x_0(4m-2)  =  m(m-4x_0)+2x_0  =  m(3n-3x_0-2) +2x_0  \hspace*{39mm}	\\
	 \hspace*{48mm}  =3m(n-x_0)-2m+2x_0= (9n-6+3x_0)(n-x_0)-6n+4. 	
	\end{eqnarray*}
	Hence, $l'_1$ must be the largest odd integer smaller than or equal to $9n -6 +3x_0 - \frac{6n-4}{n-x_0}$. 
	It is straightforward to see that $B_m(l_i')=3$ for $1\leq i\leq n-x_0$, and thus $T'$ is impossibly burnable.
	From our earlier observation, it follows that $M_n$ is bounded below by the 
	largest odd integer smaller than or equal to $12n - 2\sqrt{18n - 12} - 6$.
	
	Therefore, to complete our proof, we shall show that $l_1$ is odd for any optimal impossibly burnable $T$ with the length of its shortest path maximised. Indeed, with a more careful analysis, such $T$ would have $t_i=3$ for all $i\in B$, and furthermore,  $t_i = 4$ for all $i\in A$, assuming $n\geq 8$ for the latter as our previous observations and discussions have shown that $M_n$ is as claimed in the theorem for $3\le n \le 7$. (See \nameref{appendix B} for details.) Hence, assuming $l_1$ is not odd, it implies $l_1\geq 4m-2$.  Noting that $m\le 4n-2$,
	$$m^2-\sum_{i=1}^n l_i\le m^2- n(4m-2)=m(m-4n) +2n \leq 2n-2m<0,$$
	which is a contradiction. 
\end{proof}

\section{Conclusion}

For every $n\geq 4$, let $\Delta_n$ denote the least integer with the property that whenever $T$ is a deficient $n$-path forest with $l_1\geq \Delta_n$, then $T$ is impossibly burnable.
In our verification of Conjecture~\ref{050522f} for small values of $n$, we have observed that $\Delta_n = L_{n-1}$ for $n\in \{4,5,6,7\}$. In fact, our conjectures propose that $\Delta_n$ exists and $\Delta_n \leq L_{n-1}$ for all $n\geq 4$.
Upon further analysis using Matlab, we have found that there is only one deficient $7$-path forest with $l_1=45$ that is not impossibly burnable, namely, the path forest $(45,45,45,45,72,74,74)$, confirming $\Delta_7 = L_6=46$.
The scarceness of such deficient path forests has led us to anticipate the likelihood of $\Delta_n< L_{n-1}$ for larger $n$.
The study of the values of $\Delta_n$ potentially poses another challenging open problem.

Theorem~\ref{mainTheorem} implies that the values of $L_n$ are known if the answer to Question~\ref{mainQuestion} is affirmative.
However, although impossibly burnability is a simpler concept, Conjecture~\ref{050522f} is surprisingly nontrivial. Furthermore, as pointed out above,
$L_{n-1}$ is not necessarily the tight lower bound on the order of the shortest path for the conclusion to be true, and thus its essentiality in a possible proof by induction is doubtful.
As an alternative approach in light of Theorem~\ref{mainTheorem}, we now propose another conjecture, the truth of which implies a good asymptotic approximation to the values of $L_n$.

\begin{conjecture}\label{071222a}
	$L_n \leq 12 n$ for all $n\geq 2$. 
\end{conjecture}

Note that $L_n\geq M_n+1$ for all $n\geq 2$.
By Theorem~\ref{mainTheorem}, we have $M_n\sim 12n$.
Therefore, assuming Conjecture~\ref{071222a} holds, it follows that $L_n \sim 12n$, that is, $ \frac{L_n}{12n}\rightarrow 1$ as $n\rightarrow \infty$.

\acknowledgements
\label{sec:ack}
The second author acknowledges the support for this research by the Malaysian Ministry of Higher Education for Fundamental Research Grant Scheme with Project Code:  
FRGS/1/2023/STG06/USM/02/7.


\section*{Appendix A}\label{appendix A}

Henceforth, unless stated otherwise, $T$ is a $7$-path forest of order $m^2$ for some $m$. Note that when $l_1\geq 46$, $m$ is at least $18$. For $m$ up to $22$, we obtained the complete list of all well-burnable $7$-path forests of order $m^2$ and thus the corresponding list of deficient $7$-path forests thereafter. From here, it is easy to filter out those with $l_1\geq 46$. As a matter of fact, we saved the lists of well-burnable $7$-path forests for $m=21$ and $m=22$ in many parts, as the memory required for them to be saved as a single array is too large. Hence, we managed to verify Conjecture~\ref{050522f} (henceforth, our conjecture) for the case of seven paths for $m$ up to $22$ this way. However, we can no longer proceed in this manner for larger $m$. Therefore, in this appendix, we give a brief account on how we go around it. Table~\ref{151017b} gives some statistics from our Matlab search.

\begin{table}[h]
	\begin{center}
		\begin{tabular}{ | c | c | c| }
			\hline
			$m$ & \# well-burnable & \# deficient (impossibly burnable) \\
			&	 7-path forests with $l_1\geq 46$ & 7-path forests with $l_1\geq 46$
			\\ \hline
			18 & 2  & 0       \\ \hline
			19 & 5553 &178 \\ \hline
			20 & 162074 &1588\\ \hline
			21 & 1504741 &5460\\ \hline
			22 & 8134818 &9536\\ \hline
			23 &  31981775 & 9572\\ \hline
			24 & 101854804 &1294  \\ \hline
			25 & 279148714 & 79\\ \hline
			26 & 683537772  &  4  \\ \hline
			27 &  1532853276 &   0\\ \hline
		\end{tabular}
		\vspace{2mm}
		\caption{Verification of Conjecture~\ref{050522f} for the case of seven paths}\label{151017b}
	\end{center}
\end{table}

For convenience, we use the following definition.

\begin{definition}
	Let $n\geq 2$. 	Suppose $T =(l_1,l_2, \ldots, l_n)$ and $T'=(l'_1,l'_2, \ldots, l'_n)$ are path forests of orders $m^2$ and $(m+1)^2$, respectively.
	If there exists $1\leq i\leq n$ such that $l'_i= l_i+(2m+1)$ and $l'_j=l_j$ for all $j\neq i$, then we say that $T'$ is an \emph{extension} of $T$ (or $T$ is a \emph{reduction} of $T'$) at the $i$th component.
\end{definition}

Suppose  $\vert T\vert =23^2$ and $l_1\geq 46$. Note that $l_7\geq 76$ and thus $l_7-45\geq 31$. If $T$ is deficient, then $T'=(l_1, l_2, \ldots, l_6, l_7-45)$ is deficient. 
Hence, we first identified all  such path forests $T$ that are potentially deficient.
Such $T$ can be obtained from a deficient $7$-path forest $T'$ of order $22^2$ with $31\leq l_
1'\leq 45$ and $l_2'\geq 46$ by extension at the first component or  from a deficient $7$-path forest $T'$ of order $22^2$ with $l_1'\geq 46$ by extension at any of the seven components. This way, we obtained $36529$ potentially deficient $7$-path forests of order $23^2$ with $l_1\geq 46$.
From here, we noticed immediately that $9572$ among them are impossibly burnable. Hence, to verify our conjecture for $m=23$, it suffices to check that the remaining $26957$ path forests are all well-burnable. To show this, we first extracted from the list of all deficient $7$-path forests with $m=22$, a sublist of those with $l_1\geq 1$ and $l_2\geq 46$. We exhaustively checked and found that for each of the $26957$ path forests $T$, at least one of its seven reductions is not in the said sublist and thus $T$ is well-burnable.

To deal with larger $m$, we obtained the following two lists. 
\begin{enumerate}[{List }A.]
	\item All $9612$ deficient $7$-path forests of order $23^2$ with $l_1\geq 42$. 
	\item All $9931471$ deficient $7$-path forests of order $23^2$ with $l_1\geq 1$ and $l_2\geq 49$. 
\end{enumerate}

There are three ways to obtain a well-burnable path forest of order $23^2$ with its shortest path having order at least $42$:

\begin{enumerate}
	\item adding a path of order $45$ to any well-burnable $6$-path forest $T$ such that $\vert T\vert=22^2$ and $l_1\geq 42$;
	\item adding $45$ vertices to any of the paths of any well-burnable $7$-path forest $T$ such that $\vert T\vert=22^2$  and $l_1\geq 42$;
	\item adding $45$ vertices to the first path of any well-burnable $7$-path forest $T$ such that $\vert T\vert=22^2$, $l_1\leq 41$, and $l_2\geq 42$.
\end{enumerate}
This way, we obtained the complete list of $65485064$ well-burnable $7$-path forests of order $23^2$ with $l_1\geq 42$. From here, List A was obtained.

Before we proceed, an observation about List A will be useful later. Among the $9612$ members of List A, there are only $40$ of them with $42\leq l_1\leq 44$ and none with $l_1=45$. Furthermore, $l_7\geq 207$ for each of the $40$ path forests. (Note that $B_{23}(205)=5$ while $B_{23}(207)=7$.)

Similarly, we obtained the complete list of $311596739$ well-burnable $7$-path forests of order $23^2$ with $l_2\geq 49$. From here, List $B$ was obtained.

From both List A and List B, we extracted the sublist of deficient $7$-path forests of order $23^2$ with $l_1\geq 42$ and $l_2\geq 49$ and they coincide. Hence, this gives an assurance that our lists are correct. Furthermore, for our purposes, some sublists from the lists we obtained are probably sufficient, but as things developed, we ended up with those lists, as well as other lists that are not reported here.

Now, suppose $\vert T\vert =24^2$ and $l_1\geq 46$. Note that $l_7-47\geq 36$. If $T$ is deficient, then $T'=(l_1, l_2, \ldots, l_6, l_7-47)$ is deficient. 
Hence, similarly, we first identified all such path forests $T$ that are potentially deficient.
Such $T$ can be obtained from a deficient $7$-path forest $T'$ of order $23^2$ with $36\leq l_
1'\leq 45$ and $l_2'\geq 46$ by extension at the first component or from a deficient $7$-path forest $T'$ of order $23^2$ with $l_1'\geq 46$ by extension at any of the seven components.  
This way, we obtained $34959$ potentially deficient $7$-path forests of order $24^2$ with $l_1\geq 46$.
From here, we noticed immediately that $1294$ among them are impossibly burnable. Hence, 
to verify our conjecture for $m=24$, it suffices to check that the remaining $33665$ path forests are all well-burnable.

Note that if one of the paths of $T$ has order $47$, then deleting this path would result in a $6$-path forest that is well-burnable because $L_6=46$. Furthermore, only $34$ among those $33665$ path forests have $l_1=46$ and only one among these $34$ path forests has $l_2\neq47$, namely the path forest $(46,49,49,49,49, 92,242)$, which is \mbox{$24$-burnable} as $(46,49,49,49,49,242)$ is $22$-burnable.
Hence, filtering out those with $l_1\in \{46,47\}$, we are left with $10712$ path forests to be checked. 
If $T$ is any of the $10712$ path forests,  we observed that $l_1\geq 49$ and  found that at least one of its seven reductions is outside List B, and thus $T$ is well-burnable.

To deal with the case $m=25$, we consider the following three lists separately.
\begin{enumerate}[{List }C.]
	\item  All deficient $7$-path forests of order $25^2$ with $l_1\geq 46$ and $l_7\geq 95$. 
\end{enumerate}
\begin{enumerate}[{List }D.]
	\item All deficient $7$-path forests of order $25^2$ with $l_1\geq 46$ and $91\leq l_7\leq 94$.	
\end{enumerate}
\begin{enumerate}[{List }E.]
	\item All deficient $7$-path forests of order $25^2$ with $l_1\geq 46$ and $l_7=90$.
\end{enumerate}

Suppose $\vert T\vert=25^2$ and $l_1\geq 46$. If $l_7=90$, then $l_6=90$ and thus $T$ is well-burnable because $(l_1, \ldots, l_4, l_5-45)$ is $20$-burnable as $L_5=36$. Hence, List E is empty.

Now, suppose $91\leq l_7\leq 94$ and thus $l_6\geq 89$. Consider the path forest $T''=(l_1, \ldots, l_5, l_6-47, l_7-49)$.
If $T$ were to be deficient, then $T''$ would be among the $40$ members in List A specifically mentioned earlier as $42\leq l_1''\leq 45$. However, from our earlier observation, that would force $l_5\geq 207$, which is impossible as $l_5\leq l_7$. Hence, List D is empty as well.

To deal with List C, we first identified all $7$-path forests $T$ of order $25^2$ with $l_1\geq 46$ and $l_7\geq 95$ that are potentially deficient.
Such $T$ can be obtained from a deficient (and impossibly burnable) $7$-path forest $T'$ of order $24^2$ with $l_1'\geq 46$ by extension at any of the seven components.  This way, we obtained $5042$ path forests and $79$ among them are impossibly burnable. Hence, to verify  our conjecture for $m=25$, it suffices to check that the remaining $4963$ path forests are all well-burnable. 

Suppose $T$ is any of the $4963$ potentially deficient $7$-path forests. We noticed that the first path of $T$ is either $49,51,53,55$, or $57$. If the first path has order $49$, then deleting this path will result in a well burnable $6$-path forest  as $L_6=46$. Hence, filtering out those where $l_1=49$, we are left with $751$ path forests. If any reduction $T'$ of $T$ has $l'_1\geq  46$ but is not impossibly burnable, then $T'$ is well-burnable and thus is $T$. Filtering further based on this, we are left with a much shorter list of $45$ path forests.
Suppose $T$ is any of the $45$ path forests. We noticed that one of the paths of $T$ has order $94$. Let $T'$ be the $6$-path forest obtained from $T$ by deleting this path and deleting $47$ vertices from the longest path. Then $T'$ is well-burnable as $L_6 =46$ and thus $T$ is well-burnable. Therefore, the $79$ impossibly burnable $7$-path forests are the only deficient $7$-path forests of order $25^2$ with $l_1\geq 46$. 

Before we proceed to $m=26$, we make a note that among the $79$ impossibly burnable path forests, there are exactly four with $l_1=53$, namely:
$$(53,  53,    53,    53,    53,    53,   307),
(53,  53,    53,    53,    53,    55,   305),$$
$$(53,  53,    53,    53,    53,    57,   303),
(53,  53,    53,    53,    55,    55,   303).$$

Suppose $\vert T\vert=26^2$ and $l_1\geq 46$. If $T$ is deficient, then $T'=(l_1, \ldots, l_6, l_7-51)$ is deficient and $l_7-51\geq 46$. Hence,   we first obtained the list of $7$-path forests $T$ of order $26^2$ with $l_1\geq 46$ that are potentially deficient.
There are $292$ members in this list of potentially deficient path forests and among them, only four are impossibly burnable and they are extensions of the four above, namely:
$$(53,  53,    53,    53,    53,    53,   358),
(53,  53,    53,    53,    53,    55,   356),$$
$$(53,  53,    53,    53,    53,    57,   354),
(53,  53,    53,    53,    55,    55,   354).$$

Suppose $T$ is any of the remaining $288$ potentially deficient path forests. If the first path of $T$ has order $51$, then deleting this path will result in a well-burnable $6$-path forest  as $L_6=46$ and thus $T$ is well-burnable.
Hence, filtering out those where $l_1=51$, we are left with  $14$ path forests where $l_1=53$ for each of them coincidently. However, the reduction of each of the $14$ path forests at the last component can no longer be one of the four impossibly burnable path forests for $m=25$ with $l_1=53$. It follows that $T$ is well-burnable.

Finally, we can see that no path forest of order $27^2$ with $l_1\geq 46$ is deficient. Suppose $T$ is one such path forest and is deficient. Then $T'=(l_1, \ldots, l_6, l_7-53)$ is deficient and thus $T'$ is one of the four impossibly burnable path forests for $m=26$ above. However, this implies that $l_1=53$ (so is $l_2=l_3=53$). Since $(l_2, l_3, \ldots, l_7)$ is well-burnable as $L_6=46$, it follows that $T$ is well-burnable, which gives a contradiction. By induction, no path forest of order $m^2$  with $l_1\geq 46$ for $m\geq 27$   is deficient.

\section*{Appendix B}\label{appendix B}

In this appendix, we provide some details about the careful analysis referred to in the proof of Theorem~\ref{mainTheorem}.

Suppose there is an $i\in B$ such that $t_i\le 2$. 
Then we would have
$$m+2 \le \sum_{i=1}^n t_i \le 3(|B|-1)+2 + s + 2|A|$$ and so Inequality (\ref{mainInq}) would become
\begin{displaymath}
\left(3|B| + 2|A| - 3\right)^2 \ge l_1|B| + \left(\sum_{i,j\in A, i\ne j}s_is_j\right) + 2|A|.
\end{displaymath}
As before, it can be shown that $|B|>0$; if $|A| = 0$, then $l_1 \le 9n - 18 + \frac{9}{n}$; and if $|A| = 1$, then $l_1 \le 9n - 15 - \frac{1}{n-1}$.

If $\vert A\vert \geq 2$, letting $s= 2|A|+k$, following the same exact analysis would lead to
\begin{eqnarray*}
(3|B|+2|A|-3)^2 &\ge & l_1|B| + 4|A|^2+2k|A| - 2|A|\\
\Longrightarrow\qquad  l_1 &\le &  9n -8 +3|A| - \frac{10n-9}{n-|A|}.
\end{eqnarray*}
It is straightforward that in the range of $0<x<n$, the function $3x - \frac{10n-9}{n-x}$ is maximised when $x = n - \sqrt{\frac{10n-9}{3}}$, with the maximum value being $3n - 2\sqrt{30n-27}$.  Therefore, $l_1\le 12n - 2\sqrt{30n - 27} - 8$.

We now see that for $n\ge 3$,
\begin{eqnarray*}
l_1 & \le & \max \left\{9n-18+\frac{9}{n}, 9n-15-\frac{1}{n-1},  12n - 2\sqrt{30n - 27} - 8\right\}\\
& = & 12n - 2\sqrt{30n - 27} - 8.
\end{eqnarray*}

However, $(12n- 2\sqrt{18n - 12} - 6)-(  12n -2 \sqrt{30n-27} -8      )>4$  for $n\geq 3$.
Hence, $l_1$ is not maximised. Therefore, for an impossibly burnable $n$-path forest to be optimal such that $l_1$ is maximised, we shall need $t_i=3$ for all $i\in B$.

Now, suppose there is an $i\in A$ such that $t_i\geq 5$. We proceed from Inequality (1), namely,
\begin{displaymath}
\left(3|B| + 2|A| -2\right)^2 \ge l_1|B| + \left(\sum_{i,j\in A, i\ne j}s_is_j\right) + 2|A|.
\end{displaymath}

We have seen that $|B|>0$; if $|A| = 0$, then $l_1 \le 9n - 12 + \frac{4}{n}$; and if $|A| = 1$, then $l_1 \le 9n - 9 - \frac{2}{n-1}$. Also, if $\vert A\vert =2$, then $l_1\le 9n -6+3(2)- \frac{6n-4}{n-2}= 9n-6-\frac{8}{n-2}$.

For the final case where $|A|\ge 3$, say $s_{i_0}\geq 3$, we first observe that
\begin{eqnarray*}
\sum_{i,j\in A, i \ne j} s_is_j &= & \sum_{i\in A}s_i\left(\sum_{j\in A, j\ne i} s_j\right)
= \sum_{i\in A} s_i (s-s_i)\\
&\ge &\sum_{i\in A\backslash \{i_0\}} 2(s-2)+ s_{i_0}(s-s_{i_0}) \\
&\ge & 2(|A|-1)(s-2)+ 3(s-3) \qquad\mbox{(as $3\le s_{i_0}\le s-4$ )}\\
&= &(2 |A|+1)(s-2)-3.
\end{eqnarray*}
Letting $s = 2|A| +1+ k$ for some $k\ge 0$, we see that
\begin{displaymath}
\sum_{i,j\in A, i \ne j} s_is_j\ge (2|A|+1)(2|A|-1 +k )-3 = 4|A|^2-4+k(2|A|+1).
\end{displaymath}
Together with Inequality (\ref{mainInq}), we have
\begin{eqnarray*}
9|B|^2 + 12|A||B| - 12|B| - 10|A| + 8 &\ge  & l_1|B| + k(2|A|+1)\\
\Longrightarrow\qquad l_1 &\le &  9n -2 +3|A| - \frac{10n-8}{n-|A|}.
\end{eqnarray*}
It is straightforward that in the range of $0<x<n$, the function $3x - \frac{10n-8}{n-x}$ is maximised when $x = n - \sqrt{\frac{10n-8}{3}}$, with the maximum value being $3n - 2\sqrt{30n-24}$.  Therefore, $l_1\le 12n - 2\sqrt{30n - 24} - 2$.

We now see that for $n\ge 8$,
\begin{eqnarray*}
l_1 &\le & \max \left\{9n-12+\frac{4}{n}, 9n-9-\frac{2}{n-1}, 9n-6-\frac{8}{n-2},  12n - 2\sqrt{30n - 24} - 2\right\}\\
&= & \begin{cases}
12n - 2\sqrt{30n - 24} - 2 &\text{if } n> 8\\
9n-6 - \frac{8}{n-2} &\text{if } n=8.
\end{cases}
\end{eqnarray*}

However, $(12n- 2\sqrt{18n - 12} - 6)-(  12n -2 \sqrt{30n-24} -2      )>2$  for $n>8$ and
$(12n- 2\sqrt{18n - 12} - 6)-( 9n-6 - \frac{8}{n-2}    )>2$ for $n=8$. Hence, $l_1$ is not maximised. Therefore, for an impossibly burnable $n$-path forest to be optimal such that $l_1$ is maximised when $n\geq 8$, we shall need $t_i=4$ for all $i\in A$. (Note that $(45,45,45,45,74,107)$ is an impossibly burnable $6$-path forest such that $l_1=M_6=45$ is maximised, but $B_m(l_6)=5$ where $m=19$.)

\end{document}